\documentclass[11pt,a4paper]{amsart}

\usepackage[T1]{fontenc}
\usepackage{amssymb}
\usepackage{microtype}
\usepackage{array,booktabs,xcolor,tikz}
\usepackage[left=28mm,right=28mm,top=28mm,bottom=30mm,
            headsep=8mm,footskip=12mm]{geometry}
\usepackage[colorlinks=true,linkcolor=black,citecolor=black,urlcolor=black]{hyperref}

\makeatletter
\@namedef{subjclassname@2020}{%
  \textup{2020} Mathematics Subject Classification}

\renewcommand{\section}{%
  \@startsection{section}{1}{\z@}%
    {1.4\baselineskip \@plus .3\baselineskip \@minus .2\baselineskip}%
    {.6\baselineskip}%
    {\normalfont\large\bfseries\centering}}
\renewcommand{\subsection}{%
  \@startsection{subsection}{2}{\z@}%
    {1\baselineskip \@plus .2\baselineskip \@minus .1\baselineskip}%
    {.4\baselineskip}%
    {\normalfont\bfseries\raggedright}}

\g@addto@macro\@adjustvertspacing{%
  \abovedisplayskip=8pt plus2pt minus2pt
  \belowdisplayskip=8pt plus2pt minus2pt
  \abovedisplayshortskip=5pt plus2pt minus1pt
  \belowdisplayshortskip=5pt plus2pt minus1pt
  \jot=5pt\relax}

\makeatother

\newcommand{\Z}{\mathbb Z}
\newcommand{\R}{\mathbb R}

\newcommand{\tr}{\operatorname{tr}}
\newcommand{\diag}{\operatorname{diag}}
\newcommand{\one}{\mathbf1}
\newcommand{\eps}{\varepsilon}

\newtheorem{theorem}{Theorem}[section]
\newtheorem{lemma}[theorem]{Lemma}
\newtheorem{proposition}[theorem]{Proposition}

\theoremstyle{definition}
\newtheorem{definition}[theorem]{Definition}
\newtheorem{remark}[theorem]{Remark}
\newcommand{\reviewref}[1]{\ref{#1}}
\numberwithin{equation}{section}

\begin{document}

\title[The Rankin number \texorpdfstring{$\gamma_{5,2}$}{gamma(5,2)}]
{The Exact Value of the Rankin Number \texorpdfstring{$\gamma_{5,2}$}{gamma(5,2)}}

\author[Z. Li]{Zhixing Li}
\date{}

\begin{abstract}
The Rankin number \(\gamma_{5,2}\) measures the exact upper bound for the squared determinant of the densest $2$-dimensional sublattice of any $5$-dimensional lattice of determinant one. We prove that \(\gamma_{5,2}=16^{1/5}\), attained by $D_5^*$. 
\end{abstract}


\maketitle

\section{Preliminary}\label{sec:preliminary}

For background on lattices and positive definite quadratic forms,
see \cite[Sections~1.1 and~2.1]{Zong1999}.

\begin{definition}
Let $b_1,\ldots,b_n$ be linearly independent vectors in a Euclidean space, and write
\[
B=(b_1,\ldots,b_n),\quad
\Lambda=\Lambda(B)=B\Z^n=\Z b_1+\cdots+\Z b_n.
\]
The set $\Lambda$ is an \emph{$n$-dimensional lattice}, and the columns of $B$ form a \emph{basis} of $\Lambda$. Throughout, all vectors are column vectors.
Two basis matrices $B,B'$ generate the same lattice if and only if
\[
B'=BU,\quad U\in\mathrm{GL}_n(\Z),
\]
where $\mathrm{GL}_n(\Z)$ consists of integer matrices with determinant $\pm1$.
The \emph{determinant} of $\Lambda$ is
\[
\det(\Lambda):=\sqrt{\det(B^TB)}.
\]
Note that it is independent of the choice of the basis. If $B$ is a square matrix, then $\det(\Lambda)=|\det B|$.
The \emph{dual lattice} is
\[
\Lambda^*:=\{y\in\operatorname{span}_{\R}\Lambda:
\langle y,x\rangle\in\Z\text{ for all }x\in\Lambda\}.
\]
\end{definition}

\begin{definition}
Let
\[
\mathcal P_n:=\{P\in\mathrm{GL}_n(\R):P\text{ is positive definite and symmetric}\}.
\]
The \emph{Gram matrix} of the basis $B$ is $P=B^TB\in\mathcal P_n$. Thus
\[
\langle Bz,Bw\rangle=z^TPw\quad(z,w\in\R^n),\quad
\det P=\det(\Lambda)^2.
\]
The \emph{dual basis matrix} is $B^*=BP^{-1}$. We have
\[
(B^*)^TB=I_n,\quad \Lambda^*=B^*\Z^n,\quad
(B^*)^TB^*=P^{-1}.
\]
If $B$ is a square matrix, then $B^*=B^{-T}$.

Two lattices $\Lambda,\Lambda'$ are \emph{isometric} if $\Lambda'=O\Lambda$ for a linear Euclidean isometry $O$ between their spans. Equivalently, their basis Gram matrices satisfy
\[
P'=U^TPU,\quad U\in\mathrm{GL}_n(\Z).
\]
We call such Gram matrices \emph{arithmetically equivalent}. In particular, changing the basis from $B$ to $BU$ changes its Gram matrix from $P$ to $U^TPU$.
Two lattices are \emph{similar} if $\Lambda'=cO\Lambda$ for some $c>0$ and a linear Euclidean isometry $O$. Equivalently,
\[
P'=c^2U^TPU,\quad c>0,\quad U\in\mathrm{GL}_n(\Z).
\]
\end{definition}

We write $\mathcal L_n$ for the set of isometry classes of $n$-dimensional lattices of determinant one. An arbitrary lattice $\Lambda$ determines an element of $\mathcal L_n$ by the normalization
\[
\widehat\Lambda=(\det\Lambda)^{-1/n}\Lambda.
\]
We use the same notation for an isometry class and its representative when there is no ambiguity. Lattices used in the proof need not have determinant one; their associated points in $\mathcal L_n$ are understood through this normalization.

\begin{definition}
Let
\[
\lambda_1(\Lambda):=\min\{\|v\|:0\ne v\in\Lambda\}.
\]
The \emph{$k$-dimensional sublattice} generated by $v_1,\cdots,v_k\in\Lambda$ is $M=\Z v_1+\cdots+\Z v_k$, where the $v_i$ are linearly independent. Write $X=(v_1,\ldots,v_k)$.
\[
\det M:=\sqrt{\det(X^TX)},\quad
 d_k(\Lambda):=\min_{\substack{M\subseteq\Lambda\\\dim M=k}}\det M,
\quad 1\le k\le n.
\]
A sublattice is \emph{primitive} if 
\(M=\Lambda\cap\operatorname{span}_{\R}M\).
Every minimizing sublattice is primitive. The minimum defining \(d_k\) is attained, since basis wedges lie in the discrete exterior lattice \(\bigwedge^k\Lambda\), and their Euclidean norms are the corresponding determinants. A bounded ball contains only finitely many such wedges. See also \cite[Theorem~1, pp.~309--310]{Rankin1953}.
\end{definition}

\begin{definition}[Rankin and Hermite constants]
For \(1\le k<n\), with an abuse of notation, define the \emph{Rankin invariant} by
\[
\begin{array}{rcl}
\gamma_{n,k}(\cdot)\colon\mathcal L_n&\longrightarrow&\R\\
\Lambda &\longmapsto&
\left(\dfrac{d_k(\Lambda)}{\det(\Lambda)^{k/n}}\right)^2.
\end{array}
\]
The \emph{Rankin number} is
\(\gamma_{n,k}:=\sup\limits_{\Lambda \in\mathcal L_n}\gamma_{n,k}(\Lambda)\).
The \emph{Hermite constant} is \(\gamma_n=\gamma_{n,1}\); also write \(\gamma_n(\Lambda)=\gamma_{n,1}(\Lambda)\).
\end{definition}

\begin{remark}
Note that by the definition of the dual lattice, for \(1\le k<n\),
\[
\det(\Lambda)^{-1}=\det(\Lambda^*),\quad
 d_k(\Lambda)=d_{n-k}(\Lambda^*)\cdot\det(\Lambda).
\]
Thus \(\gamma_{n,k}=\gamma_{n,n-k}\).
\end{remark}

\begin{table}[htbp]
\caption{Known values of \(\gamma_{n,k}\) \cite{Rankin1953,Martinet2003,Sawatani2010,DutourWoerden2025,Cohn2017}.}
\label{tab:known-rankin-values}
\renewcommand{\arraystretch}{1.3}
\setlength{\tabcolsep}{7pt}
\begin{tabular}{@{}c*{9}{c}@{}}
\toprule
\(k\backslash n\)&\(2\)&\(3\)&\(4\)&\(5\)&\(6\)&\(7\)&\(8\)&\(9\)&\(24\)\\
\midrule
\(1\)&\(2/\sqrt3\)&\(2^{1/3}\)&\(\sqrt2\)&\(2^{3/5}\)&\((64/3)^{1/6}\)&\(2^{6/7}\)&\(2\)&\(2\)&\(4\)\\
\(2\)&&&\(3/2\)&&\(3^{2/3}\)&&\(3\)&&\\
\(3\)&&&&&&&\(4\)&&\\
\(4\)&&&&&&&\(4\)&&\\
\bottomrule
\end{tabular}
\end{table}

\begin{theorem}\label{thm:main}
\(\gamma_{5,2}=16^{1/5}\).
\end{theorem}

The equality in Theorem~\ref{thm:main} was conjectured in
\cite[Conjecture~6.1]{ChewOggier2026}; the lower bound was given in
\cite[Corollaries~5.8--5.9]{ChewOggier2026}. We prove the corresponding upper bound.

In the rest of this section, we give the lower bound
\(\gamma_{5,2}\ge16^{1/5}\), attained by \(D_5^*\), the dual of \(D_5\).
Section~\ref{sec:perfect-cone-reduction} classifies a local maximizer of \(\gamma_{5,2}\) into three
perfect cones associated with \(D_5,A_5,A_5^{+3}\).
Section~\ref{sec:first-cases} treats the three cases separately
and completes the proof of Theorem~\ref{thm:main}.

\begin{theorem}\label{thm:attainment}
The maximum of \(\gamma_{n,k}\) is attained by some
\(\Lambda\). Thus we can write
\[
\gamma_{n,k}=\max_\Lambda\gamma_{n,k}(\Lambda).
\]
\end{theorem}
For a proof, see \cite[Theorem~3, p.~310]{Rankin1953}.

\begin{definition}\label{def:dn}
Let \(e_i\) be the standard orthonormal vectors in $\R^n$. The \emph{root lattice} $D_n$ is defined as 
\[
D_n=\left\{(x_1,\ldots,x_n)^T\in\Z^n:
\sum x_i\equiv0\ \operatorname{mod}2\right\}
=(e_1-e_n,\ldots,e_{n-1}-e_n,2e_n)\Z^n.
\]
\(\det(D_n)=2\). Its dual lattice is \(D_n^*=\Z^n\cup\bigl(\Z^n+\dfrac12\sum e_i\bigr)\).
\end{definition}

\begin{lemma}\label{prop:d5-lower-bound}
Let
\(D_5^*=\Z^5\cup\bigl(\Z^5+(\dfrac12,\ldots,\dfrac12)^T\bigr)\).
Then \(d_2(D_5^*)=1\). Therefore, \(\gamma_{5,2}\ge\sqrt[5]{16}\). 
\end{lemma}

\begin{proof}
Let \(u,v\in D_5^*\) be linearly independent.
Let \(M=\Z u+\Z v\subset D_5^*\) be a $2$-dimensional sublattice of \(D_5^*\).
Then
\[
\det(M)^2=\|u\|^2\|v\|^2-\langle u,v\rangle^2.
\]
Note that \(\forall x=(x_1,\ldots,x_5)^T\in D_5^*\), either
\(x_i\in\Z,\ \forall i\), or \(x_i\in\Z+\dfrac12,\ \forall i\).
If \(u,v\in\Z^5\), then \(\det(M)^2\ge1\). Otherwise,
by choosing the basis \(u,u-v\) if necessary,
we may assume \(u\in\Z^5\) and \(v\in\Z^5+(\dfrac12,\ldots,\dfrac12)^T\).
Assume \(v=\dfrac12a\) with \(a_i\) odd. 
If the \(u_i\) have mixed parity, at least four minors \(u_i a_j-u_j a_i\) are odd. If they have common parity, all minors are even and at least one is nonzero. Their squared sum is at least \(4\). Thus we also have \(\det(M)^2\ge1\).

Note that \(u=(1,0,\ldots,0)^T\), \(v=(0,1,0,0,0)^T\in D_5^*\)
gives \(d_2(D_5^*)=1\).
\end{proof}

\begin{remark}
In fact, we have
\(d_k(D_n)^2=\min\{k+1,4\},\ \forall n\ge3,\ 1\le k<n\), and thus
\[
d_k(D_n^*)^2=\dfrac14d_{n-k}(D_n)^2
=\dfrac14\min\{n-k+1,4\}.
\]
\end{remark}

\section{Auxiliary inequalities and perfect-cone reduction}\label{sec:perfect-cone-reduction}

\begin{lemma}\label{lem:short-vector-bounds}
For any $5$-dimensional lattice \(\Lambda\), we have
\begin{equation}\label{eq:short-primal}
d_2(\Lambda)^2\le\sqrt2\,\lambda_1(\Lambda)^{3/2}\det(\Lambda)^{1/2},
\end{equation}
\begin{equation}\label{eq:short-dual}
d_2(\Lambda)^2\le\dfrac32\det(\Lambda)\lambda_1(\Lambda^*).
\end{equation}
\end{lemma}
\begin{proof}
Let \(v\in\Lambda\) be a shortest vector and let
\(Q=\pi_{v^\perp}(\Lambda)\) be the orthogonal projection of \(\Lambda\) onto \((\R v)^\perp\). \(Q\) is a $4$-dimensional lattice in \(v^\perp\) and \(\det\Lambda=\lambda_1(\Lambda)\cdot\det Q\). For every \(u\in\Lambda\setminus\R v\),
\[
d_2(\Lambda)\le\det(\Z u+\Z v)=\|v\|\,\|\pi_{v^\perp}(u)\|.
\]
Taking the minimum over \(u\) on the right gives
\[
\begin{aligned}
d_2(\Lambda)^2
&\le\lambda_1(\Lambda)^2\lambda_1(Q)^2\\
&\le
\sqrt2\,\lambda_1(\Lambda)^2(\det Q)^{1/2}\\
&=\sqrt2\,\lambda_1(\Lambda)^{3/2}(\det\Lambda)^{1/2}.
\end{aligned}
\]
The second inequality used the fact that $\gamma_4=\sqrt2$. 

Let \(w\in\Lambda^*\) be a shortest vector and let \(S=\Lambda\cap w^\perp\). Since
\[
\dfrac{d_2(S)^2}{\det S}\le\gamma_{4,2}=\dfrac32,
\]
we have
\[
\det\Lambda=\dfrac{\det S}{\lambda_1(\Lambda^*)}
\ge\dfrac23\dfrac{d_2(S)^2}{\lambda_1(\Lambda^*)}
\ge\dfrac23\dfrac{d_2(\Lambda)^2}{\lambda_1(\Lambda^*)}.
\]
See \cite[Theorem~6.2.1(2)]{Martinet2003} and \cite[Theorem~5]{Rankin1953} for proofs of \(\gamma_4=\sqrt2\) and \(\gamma_{4,2}=\dfrac{3}{2}\), respectively. 
\end{proof}

\begin{lemma}\label{lem:comparison}
Let $\Lambda$ be an $n$-dimensional lattice in an $n$-dimensional Euclidean
space and let $1\le k<n$. Suppose $\Lambda$ is a local maximizer of
the Rankin invariant $\gamma_{n,k}(\cdot)$. Let
$\mathcal M(\Lambda)$ be the family of densest $k$-dimensional sublattices
$M\subseteq\Lambda$, i.e. $\det(M)=d_k(\Lambda)$. Then for every invertible linear map $T$ of the ambient space, there exists
$M\in\mathcal M(\Lambda)$ such that
\begin{equation}\label{eq:comparison}
\gamma_{n,k}(\Lambda)=\dfrac{\det(M)^2}{\det(\Lambda)^{2k/n}}
\le\dfrac{\det(TM)^2}{\det(T\Lambda)^{2k/n}}.
\end{equation}
\end{lemma}
\begin{proof}
There exists \(\lambda>0\), an orthogonal matrix \(O\), and a
self-adjoint matrix \(X\) with \(\tr X=0\), such that
$T=\lambda Oe^X$. 
The two ratios in the statement remain unchanged by the scalar factor
\(\lambda\) and by \(O\). Thus we may assume
\(T=e^X\), and put \(T_t=e^{tX}\).

For \(M\in\mathcal M(\Lambda)\), define
\[
f_M(t)=\det(T_tM)^2.
\]
Then
\[
f_M(0)=\det(M)^2=d_k(\Lambda)^2,\quad
f_M(1)=\det(TM)^2.
\]
By the Cauchy--Binet formula, \(f_M\) is a nonnegative sum of
exponential functions, and hence is convex.

Suppose, on the contrary, that
\(f_M(1)<d_k(\Lambda)^2\) for every \(M\in\mathcal M(\Lambda)\).
Then for every \(t<0\),
\[
f_M(t)\ge(1-t)f_M(0)+tf_M(1)
>f_M(0)=d_k(\Lambda)^2.
\]

By finiteness of the number of sublattices with determinant below a
fixed constant, and continuity of determinant under multiplication by \(T_t\), for \(t\)
sufficiently close to zero the minimum \(d_k(T_t\Lambda)\) is attained
by some \(T_tM\) with \(M\in\mathcal M(\Lambda)\).
Consequently, for sufficiently small negative \(t\),
\[
d_k(T_t\Lambda)^2=f_M(t)>d_k(\Lambda)^2.
\]
Since \(\det(T_t\Lambda)=\det\Lambda\), this contradicts
the fact that \(\Lambda\) is a local maximizer. 
\end{proof}

Given \(P\in\mathcal P_n\), choose \(B\in\mathrm{GL}_n(\R)\) with \(P=B^TB\), and write \(\Lambda=\Lambda(B)\). Define
\[
\begin{aligned}
m(P)&:=\min\{z^TPz:0\ne z\in\Z^n\}=\lambda_1(\Lambda)^2,\\
\operatorname{Min}(P)&:=\{z\in\Z^n:z^TPz=m(P)\},\\
\mathcal V(P)&:=\left\{\sum_{z\in\operatorname{Min}(P)}c_z zz^T:c_z\ge0\right\}
\subseteq\overline{\mathcal P_n}.
\end{aligned}
\]
\[
\operatorname{Min}(\Lambda):=B\operatorname{Min}(P)
=\{v\in\Lambda:\|v\|=\lambda_1(\Lambda)\}.
\]
Note that $\operatorname{Min}(\Lambda)$ is independent of the chosen basis.

\begin{definition}
A form \(P\) is \emph{perfect} if 
$
\operatorname{span}_{\R}\{zz^T:z\in\operatorname{Min}(P)\}
=\operatorname{Sym}_n(\R).
$
A lattice \(\Lambda(B)\) is \emph{perfect} if its Gram matrix \(B^TB\) is perfect.
This definition is independent of the choice of the basis. 
\end{definition}

\begin{theorem}[Voronoi's perfect cone decomposition]\label{thm:perfect-cones}
We use the following facts without proof.
\begin{enumerate}
\item For every \(P\in\mathcal P_n\), \(\mathcal V(P)\) is a polyhedral cone.
\item For all \(P,P'\in\mathcal P_n\), \(\mathcal V(P)\cap\mathcal V(P')\) is a common face of both cones.
\item The cones of perfect forms cover \(\mathcal P_n\):
\[
\mathcal P_n\subseteq\bigcup_{P\ \mathrm{perfect}}\mathcal V(P).
\]
\item\label{item:perfect-five} Up to scaling, all $5$-dimensional perfect lattices are \(D_5,A_5,A_5^{+3}\).
\end{enumerate}
For the cone decomposition see \cite[Definition~5.1.3 and Fact~5.2.2]{ChanMeloViviani2013}. For the $5$-dimensional classification see \cite[Theorems~1.1--1.2 and p.~151]{Martinet1999}, or \cite[Corollary~6.4.3]{Martinet2003}.
\end{theorem}

We give the three perfect lattices as follows. For convenience in the following arguments, from now on, the vectors may be represented as vectors in $\R^6$, while the lattices and their corresponding quadratic forms are viewed as those restricted to the $5$-dimensional linear subspace spanned by the lattice. 
Write
\[
\one=(1,\ldots,1)^T\in\R^6,\quad E=\one^\perp,\quad
\Pi=I_6-\dfrac16\one\one^T.
\]
The restriction of \(\Pi\) to \(E\) is the identity. 
\begin{enumerate}
\item \(
D_5=\{z\in\Z^5:\textstyle\sum_i z_i\equiv0\pmod2\},\quad
B_{D_5}=(e_1-e_2,e_2-e_3,e_3-e_4,e_4-e_5,e_4+e_5),
\)
\[
\lambda_1(D_5)=\sqrt2,\quad\det(D_5)=2.
\]

\item
\(
A_5=E\cap\Z^6,\quad B_{A_5}=(e_1-e_6,\ldots,e_5-e_6),
\)
\[\lambda_1(A_5)=\sqrt2,\quad\det(A_5)=\sqrt6,
\]
\item
\(
A_5^{+3}=\sqrt{\dfrac32}\,\Pi D_6,
\quad B_{A_5^{+3}}=\sqrt{\dfrac32}\,\Pi
(e_1+e_2,e_1+e_3,e_1+e_4,e_1+e_5,e_2+e_6),
\)
\[
\lambda_1(A_5^{+3})=\sqrt2,\quad\det(A_5^{+3})=\dfrac94.
\]
Here \(D_6\) is the root lattice in Definition~\ref{def:dn}.
\end{enumerate}

Let \(\Lambda\) be a $5$-dimensional lattice. By Theorem~\reviewref{thm:perfect-cones}, one can choose a reference lattice \(\Lambda_0\in\{D_5,A_5,A_5^{+3}\}\) and bases \(B,B_0\) of $\Lambda$, $\Lambda_0$, respectively, such that
\[
P=B^TB\in\mathcal V(P_0),\text{ where } P_0=B_0^TB_0.
\]
Write
\[
P=\sum_{z\in\operatorname{Min}(P_0)}c_z zz^T,\quad c_z\ge0, 
\]
\[
\quad B_0^*=B_0P_0^{-1},\quad\Lambda_0^*=B_0^*\Z^5,
\]
and put \(K=B_0PB_0^T\).
Still, all powers, inverses, traces, and determinants of \(K\) are taken on the $5$-dimensional Euclidean span of \(\Lambda_0\). Up to an orthogonal isometry, we may realize
\[
\Lambda=K^{1/2}\Lambda_0^*,\quad B=K^{1/2}B_0^*.
\]
\begin{proposition}\label{lem:reference-bounds}
\begin{equation}\label{eq:ref-det}
\det K=(\det\Lambda)^2(\det\Lambda_0)^2,
\end{equation}
\begin{equation}\label{eq:ref-trace}
\tr K\le\lambda_1(\Lambda_0)^2\cdot\dfrac{45}{4}\cdot
\dfrac{\det(\Lambda)^2}{d_2(\Lambda)^4}.
\end{equation}
\end{proposition}
\begin{proof}
Equation \eqref{eq:ref-det} follows from $K=B_0B^TBB_0^T$. 
Since
\(
K=\sum_{z\in\operatorname{Min}(P_0)}c_z(B_0z)(B_0z)^T, 
\)
we have
\[
\tr K=\lambda_1(\Lambda_0)^2\sum_zc_z.
\]
Since \(P^{-1}\) is the Gram matrix of the basis of \(\Lambda^*\) dual to \(B\),
\[
\begin{aligned}
5=\tr(PP^{-1})
&=\sum_zc_z\,z^TP^{-1}z\\
&\ge\lambda_1(\Lambda^*)^2\sum_zc_z\\
&\overset{\eqref{eq:short-dual}}{\ge}
\dfrac49\dfrac{d_2(\Lambda)^4}{\det(\Lambda)^2}\sum_zc_z.
\end{aligned}
\]
Combining the last two displays proves \eqref{eq:ref-trace}.
\end{proof}

\section{Proof of the main theorem}\label{sec:first-cases}
\begin{proof}[Proof of Theorem~\reviewref{thm:main}]
We prove the theorem by contradiction. Suppose that
\(\gamma_{5,2}>\sqrt[5]{16}\), and, by Theorem~\reviewref{thm:attainment}, we may choose a maximizer \(\Lambda\) of $\gamma_{5,2}(\cdot)$.
Scale \(\Lambda\) so that \(d_2(\Lambda)=1\). Then
$\det\Lambda<\dfrac12$. 
From now on, let \(\Lambda_0,B_0,P_0,B_0^*,P,B,K\) be defined as in the previous section. 
Equations \eqref{eq:short-primal} and \eqref{eq:short-dual} imply
\begin{equation}\label{eq:floors}
\lambda_1(\Lambda)>1,\quad
\det(\Lambda)\lambda_1(\Lambda^*)\ge\dfrac23.
\end{equation}

\subsection*{\texorpdfstring{Case 1. $\Lambda_0=A_5$}{Case 1. Lambda0 = A5}}

Denote
\[
d_i:=\Pi e_i=e_i-\dfrac16\one,\quad 1\le i\le6.
\]
\(A_5^*=\Z d_1+\cdots+\Z d_6\).
Write \(b_i=K^{1/2}d_i\in E\). Thus \(\Lambda=\Z b_1+\cdots+\Z b_6\), and any five of these vectors form a basis.
For \(i\ne j\), \(\Z b_i+\Z b_j\) is a $2$-dimensional sublattice of \(\Lambda\), and
\[
\det(\Z b_i+\Z b_j)\ge d_2(\Lambda)=1.
\]
Note that
\[
\sum_{i=1}^6d_id_i^T=\Pi,\quad
\sum_{i=1}^6b_ib_i^T=K.
\]
Consequently, by the Cauchy--Schwarz inequality, 
\[
\begin{aligned}
15
&\le\sum_{i<j}\|b_i\wedge b_j\|^2
=\sum_{i<j}\bigl(\|b_i\|^2\|b_j\|^2-\langle b_i,b_j\rangle^2\bigr)\\
&=\dfrac12\bigl((\tr K)^2-\tr(K^2)\bigr)
\le\dfrac25(\tr K)^2
\overset{\eqref{eq:ref-trace}}{<}\dfrac{405}{32}<15.
\end{aligned}
\]
This is a contradiction; thus \(\Lambda_0\ne A_5\).


\subsection*{\texorpdfstring{Case 2. $\Lambda_0=A_5^{+3}$}{Case 2. Lambda0 = A5+3}}

Recall that
\[
\Lambda_0=A_5^{+3}=\sqrt{\dfrac32}\,\Pi D_6,
\text{ where }
D_6:=\{z\in\Z^6:\sum_i z_i\equiv0\pmod2\},
\]
and take
\[
B_0:=B_{A_5^{+3}}
=\sqrt{\dfrac32}\,\Pi(e_1+e_2,e_1+e_3,e_1+e_4,e_1+e_5,e_2+e_6),
\quad P_0:=B_0^TB_0.
\]
It is easy to verify the following (for the classical Coxeter construction and shortest-vector description, see \cite[Proposition~5.2.3]{Martinet2003}): 

\begin{enumerate}
\item
\(
\Lambda_0^*=\sqrt{\dfrac23}
\{v\in\Z^6\cup(\Z+\dfrac12)^6:\sum_i v_i=0\},
\quad \lambda_1(\Lambda_0^*)=1;
\)
\item $\lambda_1(\Lambda_0)=\sqrt2, \quad \det(\Lambda_0)=\dfrac{9}{4}$.
\item 
\(
\operatorname{Min}(\Lambda_0)
=\{\pm\sqrt{\dfrac{3}{2}}\,\Pi(e_i+e_j):1\le i<j\le6\},
\quad \#\operatorname{Min}(\Lambda_0)=30.
\)
\item $\operatorname{Min}(\Lambda_0^*)$ consists of all coordinate
permutations of
\[
\dfrac1{\sqrt6}(1,1,1,-1,-1,-1)^T,
\quad \#\operatorname{Min}(\Lambda_0^*)=20.
\]
\end{enumerate}

Recall that the basis $B$ of $\Lambda$ satisfies 
\[
P=B^TB\in\mathcal V(P_0), \quad
P=\sum_{z\in\operatorname{Min}(P_0)}c_z\,zz^T,
\ c_z\ge0, \quad
\text{and }K=B_0PB_0^T. 
\]

By \eqref{eq:ref-det} and
\eqref{eq:ref-trace},
\[
\det K<\dfrac{81}{64},\quad \tr K<\dfrac{45}{8}.
\]

\begin{proposition}\label{prop:a5-short}
Define
\[
\begin{aligned}
s_1&=\dfrac1{\sqrt6}(1,1,1,-1,-1,-1)^T,\\
s_2&=\dfrac1{\sqrt6}(2,0,0,0,0,-2)^T,\\
s_3&=\dfrac1{\sqrt6}(3,1,-1,-1,-1,-1)^T.
\end{aligned}
\]
If $v\in\Lambda_0^*\setminus\{0\}$ and
$\|K^{1/2}v\|^2<2$, then $v$ is of type $s_i$, $i=1,2,3$. 
Here and below, a vector is said to be \emph{of type $s_i$} if it is obtained from $s_i$ by a coordinate permutation and an overall change of sign. 
\end{proposition}

\begin{proof}
The map $v\mapsto K^{1/2}v$ is a one-to-one correspondence from
$\Lambda_0^*$ to $\Lambda$. For every $v\in\Lambda_0^*$,
\[
\begin{aligned}
\|K^{1/2}v\|^2
&=v^TKv=v^TB_0PB_0^Tv\\
&=v^TB_0\left(\sum_{z\in\operatorname{Min}(P_0)}c_z\,zz^T\right)B_0^Tv\\
&=v^T\left(\sum_{x\in\operatorname{Min}(\Lambda_0)}c_x\,xx^T\right)v\\
&=\sum_{x\in\operatorname{Min}(\Lambda_0)}c_x(x^Tv)^2\\
&=\dfrac32\sum_{i<j}c_{ij}
\bigl((e_i+e_j)^T\Pi v\bigr)^2\\
&=\dfrac32\sum_{i<j}c_{ij}
\left(v_i+v_j-\dfrac13\sum_kv_k\right)^2. 
\end{aligned}
\]
Here $x=B_0z$, $c_x=c_z$, and
\[
x_{ij}=\sqrt{\dfrac{3}{2}}\,\Pi(e_i+e_j),\quad
c_{ij}:=c_{x_{ij}}+c_{-x_{ij}}.
\]
Since $v\in E$, $\sum_kv_k=0$. Write
$v=\dfrac{z}{\sqrt6}$, where $z$ has integer coordinates
of the same parity and $\sum_i z_i=0$. Then
\begin{equation}\label{eq:a5third-energy}
\|K^{1/2}(\dfrac{z}{\sqrt6})\|^2
=\dfrac14\sum_{i<j}c_{ij}(z_i+z_j)^2.
\end{equation}

We choose $w\in\sqrt6\,\Lambda_0^*$ such that, for every $i<j$,
\[
\min\{0,z_i+z_j\}\le w_i+w_j\le\max\{0,z_i+z_j\}.
\]
By a permutation of coordinates, we may assume that $z=(a,b,c,d,e,f)^T$ with $a\ge b\ge c\ge d\ge e\ge f$. Also allowing an overall change of sign and reordering, it is easy to check that we have only three possibilities for $z$ and we assign a $w$ for each case: 
\begin{enumerate}
\item $b+c>0$ and $d+e<0$. Take $w=(1,1,1,-1,-1,-1)^T$.
\item $b+c\le0$ and $b+f<0$. Take $w=(2,0,0,0,0,-2)^T$. 
\item $b+c\le0$ and $b+f=0$. Take $w=(3,1,-1,-1,-1,-1)^T$.
\end{enumerate}
Now we have
\[
\begin{aligned}
2>\|K^{1/2}(\dfrac{z}{\sqrt6})\|^2
&=\dfrac14\sum_{i<j}c_{ij}(z_i+z_j)^2\\
&\ge\dfrac14\sum_{i<j}c_{ij}(w_i+w_j)^2
 +\dfrac14\sum_{i<j}c_{ij}\bigl((z-w)_i+(z-w)_j\bigr)^2\\
&=\|K^{1/2}(\dfrac{w}{\sqrt6})\|^2
 +\|K^{1/2}(\dfrac{z-w}{\sqrt6})\|^2.
\end{aligned}
\]
Both $\dfrac{w}{\sqrt6}$ and $\dfrac{z-w}{\sqrt6}$ belong to $\Lambda_0^*$.
Therefore $z=w$, since otherwise $\lambda_1(\Lambda)>1$ would force 
the right side to be greater than $2$. Thus $v$ is
of one of the three stated types. 
\end{proof}

\begin{proposition}\label{clm:a5-gauss}
Let $M\subset\Lambda$ be a $2$-dimensional densest sublattice, and $\{K^{1/2}u,K^{1/2}v\}$ be a Lagrange reduced basis of $M$. 
Here a basis $\{p,q\}$ of a $2$-dimensional lattice is \emph{Lagrange reduced} if
$\|p\|\le\|q\|$ and $2|\langle p,q\rangle|\le\|p\|^2$.
Then $u$ is of type $s_1$, and $v$ is of type $s_1$ or $s_2$.
\end{proposition}

\begin{proof}
Write
\[
a=\|K^{1/2}u\|^2,\quad b=\|K^{1/2}v\|^2,\quad t=u^TKv.
\]
Since $\det M=1$ and the basis is Lagrange reduced, we have
\[
1<a\le b,\quad 2|t|\le a,\quad ab-t^2=1.
\]
Thus $a^2\le ab\le1+\dfrac{a^2}{4}$, and
\[
a\le\dfrac2{\sqrt3}<\dfrac76,\quad
b=\dfrac{1+t^2}{a}\le\dfrac1a+\dfrac a4<\dfrac54.
\]
By Proposition~\ref{prop:a5-short}, both $u$ and $v$ are of type
$s_1,s_2$, or $s_3$. We prove that
\begin{equation}\label{eq:a5third-family-bounds}
\begin{aligned}
\|K^{1/2}r\|^2&>\dfrac76 &&(r\text{ of type }s_2),\\
\|K^{1/2}h\|^2&>\dfrac{73}{48}>\dfrac32
&&(h\text{ of type }s_3).
\end{aligned}
\end{equation}
All coordinate permutations below also relabel the coefficients $c_{ij}$.

First let $r$ be of type $s_2$. By a permutation and an overall
change of sign, write
\[
p=e_1-e_2,\quad r=\dfrac{2p}{\sqrt6},\quad
\rho=\|K^{1/2}r\|^2.
\]
Let
\[
\begin{aligned}
t_1&=e_3+e_4-e_5-e_6,\\
t_2&=e_3+e_5-e_4-e_6,\\
t_3&=e_3+e_6-e_4-e_5.
\end{aligned}
\]
For each $j=1,2,3$, the two vectors
$\dfrac{p+t_j}{\sqrt6},\dfrac{-p+t_j}{\sqrt6}$ belong to $\Lambda_0^*$.
Since $d_2(\Lambda)=1$,
\[
\begin{aligned}
1&\le\det\left(\Z K^{1/2}\dfrac{p+t_j}{\sqrt6}
 +\Z K^{1/2}\dfrac{-p+t_j}{\sqrt6}\right)^2\\
&=\dfrac19\bigl((p^TKp)(t_j^TKt_j)-(p^TKt_j)^2\bigr).
\end{aligned}
\]
Note that $p^TKp=\dfrac{3\rho}{2}$. Thus
$t_j^TKt_j\ge\dfrac6\rho$. Using
$K=\sum_{k<l}c_{kl}x_{kl}x_{kl}^T$, we have
\[
\rho=\sum_{k=3}^6(c_{1k}+c_{2k}),\quad
\tr K=2c_{12}+\dfrac\rho2+\dfrac13\sum_{j=1}^3t_j^TKt_j.
\]
Here $\sum_{j=1}^3(t_j^Tx_{kl})^2$ equals $0$, $\dfrac92$, or $6$,
according as both indices, one index, or neither index belongs to
$\{1,2\}$. Consequently,
\[
\tr K\ge\dfrac\rho2+\dfrac6\rho.
\]
If $\rho\le\dfrac76$, the right side is at least
$\dfrac{481}{84}>\dfrac{45}{8}$, a contradiction. This proves the
first inequality in \eqref{eq:a5third-family-bounds}.

Next let $h$ be of type $s_3$. By a permutation and an overall
change of sign, we may assume that $h=s_3$. For $k=3,4,5,6$, put
\[
u_k=\dfrac{2(e_1+e_2+e_k)-\one}{\sqrt6},\quad
r_k=\dfrac{2(e_1-e_k)}{\sqrt6}.
\]
Then $u_k$ is of type $s_1$, and $r_k$ is of type $s_2$.
By \eqref{eq:a5third-energy},
\[
\|K^{1/2}h\|^2
=\dfrac12\sum_{k=3}^6
 \left(\|K^{1/2}u_k\|^2+\|K^{1/2}r_k\|^2\right)
-\dfrac12\tr K+c_{12}.
\]
Indeed, the coefficients of $c_{12},c_{1k},c_{2k}$ and $c_{kl}$
with $k,l\ge3$ are respectively $4,1,0,1$ on both sides.
Since $\lambda_1(\Lambda)>1$ and the first inequality in
\eqref{eq:a5third-family-bounds} gives
$\|K^{1/2}r_k\|^2>\dfrac76$, we obtain
\[
\|K^{1/2}h\|^2>
\dfrac12\left(4+4\cdot\dfrac76\right)-\dfrac{45}{16}
=\dfrac{73}{48}>\dfrac32.
\]
This proves \eqref{eq:a5third-family-bounds}.

The bound $a<\dfrac76$ excludes type $s_2$ for $u$, and
$a,b<\dfrac54<\dfrac32$ excludes type $s_3$ for both $u$ and $v$.
Thus $u$ is of type $s_1$, and $v$ is of type $s_1$ or $s_2$.
\end{proof}

Let $M\subset\Lambda$ be any densest $2$-dimensional sublattice,
and choose a Lagrange reduced basis $\{K^{1/2}u,K^{1/2}v\}$.
By Proposition~\ref{clm:a5-gauss}, $u$ is of type $s_1$, and $v$ is
of type $s_1$ or $s_2$. We show that $M$ has a basis of the form
$\{K^{1/2}u,K^{1/2}v\}$ with both $u,v$ of type $s_1$.
If $v$ is of type $s_1$, there is nothing to prove. Otherwise,
by a permutation and an overall change of sign, write
\[
u=\dfrac{(1,1,1,-1,-1,-1)^T}{\sqrt6},\quad
v=\dfrac{2(e_i-e_j)}{\sqrt6}.
\]
If $i,j$ belong to different triples $\{1,2,3\}$ and $\{4,5,6\}$,
one of $u+v,u-v$ is of type $s_1$. Replacing $v$ by this vector
preserves the lattice $\Z u+\Z v$, and gives the required basis.

Suppose that $i,j$ belong to the same triple. Then both $u+v,u-v$
are of type $s_3$. With $a,b,t$ as in the preceding proof, we have
\[
1<a\le\dfrac2{\sqrt3},\quad \dfrac76<b<\dfrac54,\quad ab-t^2=1.
\]
Therefore
\[
\begin{aligned}
\min\{\|K^{1/2}(u+v)\|^2,\|K^{1/2}(u-v)\|^2\}
&=a+b-2|t|\\
&=a+b-2\sqrt{ab-1}\\
&\le\dfrac{13}{6}-\dfrac2{\sqrt6}<\dfrac32.
\end{aligned}
\]
For the estimate, $a+b-2\sqrt{ab-1}$ is decreasing in each
variable on the stated intervals: its partial derivatives are
$1-\dfrac b{\sqrt{ab-1}}<0$ and
$1-\dfrac a{\sqrt{ab-1}}<0$, since $\sqrt{ab-1}<\dfrac34<a,b$.
Its value at $a=1,b=\dfrac76$ is
$\dfrac{13}{6}-\dfrac2{\sqrt6}$.
This contradicts \eqref{eq:a5third-family-bounds}.

Thus every densest $2$-dimensional sublattice $M\subset\Lambda$ can
be written as
\[
M=\Z(K^{1/2}u)+\Z(K^{1/2}v)=K^{1/2}(\Z u+\Z v),\quad
u,v\in\operatorname{Min}(\Lambda_0^*).
\]
Note that for linearly independent vectors $u,v$ of type $s_1$,
$\|u\|=\|v\|=1$ and $|\langle u,v\rangle|=\dfrac13$. Thus
\[
\det(K^{-1/2}M)^2=1-\dfrac19=\dfrac89.
\]
Taking $T=K^{-1/2}$ in Lemma~\ref{lem:comparison}, there exists a
densest $2$-dimensional sublattice $M\subset\Lambda$ such that
\[
\sqrt[5]{16}<\gamma_{5,2}
\le\dfrac{\det(K^{-1/2}M)^2}{\det(\Lambda_0^*)^{4/5}}
=\left(\dfrac94\right)^{4/5}\dfrac89
=\left(\dfrac{128}{9}\right)^{1/5}<\sqrt[5]{16}.
\]
This is a contradiction; thus $\Lambda_0\ne A_5^{+3}$.

\subsection*{\texorpdfstring{Case 3. $\Lambda_0=D_5$}{Case 3. Lambda0 = D5}}

The argument is similar to the preceding case. Proposition~\ref{clm:d5-coverage}
classifies the possible sublattices $K^{-1/2}M$ for a densest
$2$-dimensional sublattice $M\subset\Lambda$, and shows that one of them
is of signed $2+3$ type. For such an $M$, Proposition~\ref{clm:d5-phase}
gives bounds on $\tr K$ and on the determinant after orthogonal projection.
Proposition~\ref{clm:d5-dual} then gives an inequality involving a dual
vector, which yields a contradiction when combined with these bounds.
The proof of Proposition~\ref{clm:d5-coverage} is given in Appendix~\ref{app:d5};
the proofs of Propositions~\ref{clm:d5-phase} and \ref{clm:d5-dual}
are given in Appendix~\ref{app:phase}.

Take $B_0=B_{D_5}$ as in Section~\ref{sec:perfect-cone-reduction},
and put $P_0=B_0^TB_0$. We have
\[
\Lambda_0^*=D_5^*=\Z^5\cup(\Z+\dfrac12)^5,\quad
B_0^*=B_0^{-T},\quad \det(\Lambda_0^*)=\dfrac12.
\]
Recall that
\[
d_2(\Lambda)=1,\quad \det\Lambda<\dfrac12,\quad
\lambda_1(\Lambda)>1,
\]
and
\[
P=B^TB\in\mathcal V(P_0),\quad K=B_0PB_0^T,\quad
\Lambda=K^{1/2}\Lambda_0^*,\quad \Lambda^*=K^{-1/2}\Lambda_0.
\]
Let
\[
\mathcal R=\{e_i-e_j,e_i+e_j:1\le i<j\le5\}.
\]
For $r\in\mathcal R$, write $z_r=B_0^{-1}r$ and
$c_r=c_{z_r}+c_{-z_r}$. The cone expansion gives
\[
K=\sum_{r\in\mathcal R}c_r rr^T,\quad c_r\ge0,\quad
\tr K=2\sum_{r\in\mathcal R}c_r.
\]
Thus $K$ is \emph{symmetrically diagonally dominant} (SDD).
Equations~\eqref{eq:ref-det} and \eqref{eq:ref-trace} give
\[
\det K=4\det(\Lambda)^2<1,\quad
\tr K\le\dfrac{45}{8}\det K<\dfrac{45}{8}.
\]

For a nonempty proper subset $A\subset\{1,\ldots,5\}$ and
$\sigma_i\in\{\pm1\}$, write
\[
u_A=\sum_{i\in A}\sigma_i e_i,\quad
u_{A^c}=\sum_{i\notin A}\sigma_i e_i,\quad
N_{A,\sigma}=\Z\dfrac{u_A+u_{A^c}}2+
                 \Z\dfrac{u_A-u_{A^c}}2.
\]
These are $2$-dimensional sublattices of $\Lambda_0^*$.

\begin{proposition}\label{clm:d5-coverage}\label{clm:d5-contact}
Let $M\subset\Lambda$ be a densest $2$-dimensional sublattice.
Then $K^{-1/2}M$ is one of the following $130$ primitive $2$-dimensional
sublattices of $\Lambda_0^*$. Their squared determinants are given in the last column.
\[
\begin{array}{c|c|c|c}
\text{type}&\text{reference sublattice}&\text{number}&\det^2\\\hline
\text{coordinate}&\Z e_i+\Z e_j\ (i<j)&10&1\\
\text{signed }1+4&N_{A,\sigma},\ |A|=1&40&1\\
\text{signed }2+3&N_{A,\sigma},\ |A|=2&80&\dfrac32
\end{array}
\]
Moreover, there exists a densest sublattice $M$ such that
$K^{-1/2}M=N_{A,\sigma}$ with $|A|=2$.
After a signed permutation of coordinates, we may write
\[
M=K^{1/2}M_0,\quad M_0=\Z h_++\Z h_-,
\]
where
\[
u_A=e_1+e_2,\quad u_B=e_3+e_4+e_5,\quad
h_+=\dfrac{u_A+u_B}{2},\quad h_-=\dfrac{u_A-u_B}{2}.
\]
\end{proposition}

Fix $M,M_0$ as in Proposition~\ref{clm:d5-coverage}, and write
\[
m_\pm=K^{1/2}h_\pm,\quad x_i=K^{1/2}e_i,\quad
y_i=\pi_{(\operatorname{span}_{\R}M)^\perp}(x_i).
\]
Since $x_1+x_2=m_++m_-$ and $x_3+x_4+x_5=m_+-m_-$,
\[
y_1+y_2=0,\quad y_3+y_4+y_5=0.
\]
Put
\[
A_0=\|x_1+x_2\|^2,\quad
A_B^2=\det(\Z y_3+\Z y_4)^2,\quad
s=\dfrac{45}{8}-\tr K>0.
\]

\begin{samepage}
\begin{proposition}\label{clm:d5-phase}
With the notation above,
\begin{equation}\label{eq:phaseinputs}
\tr K\ge\dfrac{71\sqrt2}{18},\quad A_0<\dfrac32,
\end{equation}
and
\begin{equation}\label{eq:d5-area-gap}
2A_B^2-\dfrac89>\dfrac{71}{24336}+\dfrac52s>0.
\end{equation}
\end{proposition}
\end{samepage}

\begin{proposition}\label{clm:d5-dual}
Let $c_-=c_{e_1-e_2}$. Then
\begin{equation}\label{eq:rootmass}
2c_->\dfrac7{12},
\end{equation}
and
\begin{equation}\label{eq:new-d5-trace}
5\det K\ge\dfrac89\tr K+2c_-\left(2A_B^2-\dfrac89\right).
\end{equation}
\end{proposition}

Combining these two propositions and using $\tr K=\dfrac{45}{8}-s$, we obtain
\[
\begin{aligned}
5\det K
&>5-\dfrac89s+
\dfrac7{12}\left(\dfrac{71}{24336}+\dfrac52s\right)\\
&=5+\dfrac{497}{292032}+\dfrac{41}{72}s>5.
\end{aligned}
\]
This is a contradiction; thus $\Lambda_0\ne D_5$.
Hence $\gamma_{5,2}\le16^{1/5}$. Together with
Lemma~\ref{prop:d5-lower-bound}, we obtain
\[
\gamma_{5,2}=16^{1/5},\quad
d_2(\Lambda)^5\le4\det(\Lambda)^2
\]
for every $5$-dimensional Euclidean lattice.
\end{proof}

\begin{remark}
Proposition~\ref{clm:d5-dual} uses the following identities.
For a positive definite $n\times n$ matrix $K=\sum_i c_i r_i r_i^T$, $c_i\ge0$,
\[
\sum_i c_i\|r_i\|^2=\tr K,\quad
\sum_i c_i r_i^TK^{-1}r_i=n.
\]
In the proof above, we use a lower bound on $r^TK^{-1}r$ for every root
$r$, and a stronger bound for $r=e_1-e_2$. The latter follows from the
determinant of a primitive sublattice.
\end{remark}

\appendix

\section{Classification of densest sublattices}\label{app:d5}

\begin{proof}[Proof of Proposition~\ref{clm:d5-coverage}]
We use the notation and assumptions of the $D_5$ case.
For $x\in\R^5$, write $q(x)=\|K^{1/2}x\|^2$.
Since $K$ is SDD, we have
\begin{equation}\label{eq:sdd}
q(x)=\sum_i d'_i x_i^2+\sum_{i<j}|K_{ij}|
(x_i+\operatorname{sgn}(K_{ij})x_j)^2,\quad d'_i\ge0.
\end{equation}
For $a>0$, let $\phi_a(t)=\max(-a,\min(t,a))$ and
$\psi_a(t)=t-\phi_a(t)$. Both functions are nondecreasing and odd;
we apply them coordinatewise to vectors.
Thus $(\phi_a(t)\pm\phi_a(u))(\psi_a(t)\pm\psi_a(u))\ge0$.
By \eqref{eq:sdd}, we have
\begin{equation}\label{eq:firm}
q(x)\ge q(\phi_a(x))+q(x-\phi_a(x)).
\end{equation}
Put $d_i=K_{ii}$ and $\nu=\lambda_1(\Lambda)^2>1$.
The images of the ten coordinate sections and the forty signed
$1+4$ sections have determinant at least one, since $d_2(\Lambda)=1$.
Also $\tr K<\dfrac{45}{8}$. Every proper signed indicator
$w=\sum_{j\in A}\sigma_j e_j$ of support size $2,3$, or $4$ satisfies
\begin{equation}\label{eq:indicator}
q(w)>\dfrac54.
\end{equation}
Indeed, put $C=A^c$, $W=q(w)$, $a=\sum_A d_i$, and $c=\sum_C d_i$.
For $|A|\ge3$, average the full half-sign vectors with the signs
on $A$ fixed. This gives $W+c\ge4\nu$, while $a\ge|A|\nu$, so
\[
W\ge(|A|+4)\nu-\tr K>\dfrac{11}{8}>\dfrac{5}{4}.
\]
For $A=\{1,2\}$, suppose $W\le\dfrac{5}{4}$. The coordinate-plane
inequality gives
\[
1\le d_1d_2-K_{12}^2\le \dfrac{aW}{2}-\dfrac{W^2}{4},
\quad a\ge\dfrac2W+\dfrac W2\ge\dfrac{89}{40},\quad c\le\dfrac{17}{5}.
\]
For each $i\in C$, average the signed $1+4$ plane inequalities through $e_i$,
fixing the signs on $A$ and averaging the other two signs.
Dropping their nonnegative squares gives $d_i(W+c-d_i)\ge4$.
Summation and $\sum_Cd_i^2\ge \dfrac{c^2}{3}$ now yield
\[
12\le cW+\dfrac23c^2
\le\dfrac54\dfrac{17}{5}+\dfrac23\left(\dfrac{17}{5}\right)^2
=\dfrac{3587}{300}<12.
\]

An integral vector which is not a
signed indicator clips at one into two nonzero integral vectors;
\eqref{eq:firm} gives norm squared at least $2\nu>2$.
A half-integral vector which is not a full half-sign vector clips at
$\dfrac{1}{2}$ into a full half-sign vector and a nonzero integral residual,
with the same lower bound. Proper non-axis indicators are excluded
below $\dfrac{5}{4}$ by \eqref{eq:indicator}. A full integer indicator is twice
a half-sign vector and has squared norm greater than four.
Consequently every nonzero $x\in D_5^*$ satisfies
\[
q(x)<\dfrac{5}{4}\ \Longrightarrow\ x=\pm e_i\quad\text{or}\quad
x=\dfrac{(\sigma_1,\ldots,\sigma_5)^T}{2}.
\]
Let $\{K^{1/2}u,K^{1/2}v\}$ be a Lagrange reduced basis of a
$2$-dimensional sublattice with determinant at most one, and write
$a=\|K^{1/2}u\|^2$, $b=\|K^{1/2}v\|^2$, and $t=u^TKv$.
Then
\[
1<a\le b,\quad 2|t|\le a,\quad ab-t^2\le1,
\]
and hence
\[
a\le\dfrac2{\sqrt3},\quad
b\le\dfrac1a+\dfrac a4<\dfrac54.
\]
Thus $u,v$ are two signed coordinate vectors, a signed coordinate vector
and a full half-sign vector, or two full half-sign vectors. The first two cases give coordinate and signed $1+4$ planes.
In the last case, let $A$ consist of the coordinates where the
two half-sign vectors agree. With $u_A,u_{A^c}$ as above, the basis is
$\dfrac{u_A+u_{A^c}}2,\dfrac{u_A-u_{A^c}}2$, with wedge
$-\dfrac{u_A\wedge u_{A^c}}2$.
The section is primitive: both block coefficients are simultaneously
integral or simultaneously half-integral, exactly as generated by this basis. The block sizes are $1+4$ or $2+3$. This proves the
classification of all $2$-dimensional sublattices of determinant at most one. For a block of size $k=1,2$,
there are $\binom5k$ choices of support and $\dfrac{2^5}{4}=8$ distinct sign
patterns modulo independent reversals of the two blocks. Thus the
numbers are $40$ and $80$. Their reference squared determinants are
$\dfrac{k(5-k)}{4}$, namely $1$ and $\dfrac{3}{2}$; the ten coordinate sections have
squared determinant one. This proves the classification assertion
of Proposition~\reviewref{clm:d5-coverage}.

It remains to show that $K^{-1/2}M$ is of signed $2+3$ type for some
densest sublattice $M$. Suppose, on the contrary, that every densest
sublattice $M$ has $K^{-1/2}M$ of one of the first two types in the table. Taking $T=K^{-1/2}$ in
Lemma~\ref{lem:comparison}, we obtain
\[
\det(\Lambda)^{-4/5}\le\det(\Lambda_0^*)^{-4/5},
\quad\text{hence}\quad (\det K)^{-2/5}\le1,
\]
since all these reference sections have squared determinant one.
This contradicts $\det K<1$. Thus $K^{-1/2}M$ is of signed $2+3$ type for some densest sublattice $M$. A signed permutation of coordinates gives
the stated $M_0$. It also conjugates $K$ and relabels the coefficients
$c_r$, preserving the root expansion and the trace and determinant bounds.

\end{proof}

\section{\texorpdfstring{The signed $2+3$ case}{The signed 2+3 case}}\label{app:phase}

\begin{proof}[Proof of Proposition~\ref{clm:d5-phase}]
We use the notation of the $D_5$ case. We first give two expressions
for $\tr K$. Let
\[
G_M=\begin{pmatrix}a&c\\c&b\end{pmatrix}
=\operatorname{Gram}(m_+,m_-),\quad ab-c^2=1.
\]
A Lagrange reduced basis of $M$ has both squared norms below $\dfrac{5}{4}$. The reference plane of $M_0$ contains no coordinate axis, and its
only full half-sign lines are $\R h_+$ and $\R h_-$.
By Appendix~\ref{app:d5}, every vector in $M$ of squared norm below
$\dfrac54$ lies on $\R m_+$ or $\R m_-$. Changing the signs of
coordinates $3,4,5$ interchanges $h_+$ and $h_-$, and conjugates $K$ by
the same signed permutation.
We may therefore arrange
\[
1<a\le b,\quad |c|\le \dfrac a2,\quad a\le\dfrac{2}{\sqrt3},\quad b<\dfrac{5}{4},
\quad |c|\le\dfrac{1}{\sqrt3}.
\]
Write $x_i=r_i m_++s_i m_-+y_i$, with $y_i\perp\operatorname{span}_{\R}M$, and put
$z_i=(r_i,s_i)^T$, $q_i=\|y_i\|^2$. The coordinate relations imply
\[
\sum_{i=1}^2z_i=(1,1)^T,\quad
\sum_{i=3}^5z_i=(1,-1)^T.
\]
The sublattices $\Z x_i+\Z m_\pm\subseteq\Lambda$, $i=1,\ldots,5$,
have squared determinants
\[
F_i^+=a q_i+s_i^2\ge1,\quad F_i^-=b q_i+r_i^2\ge1.
\]
Set
\[
R=a+b,\quad \ell=\dfrac1{2a}+\dfrac1{2b},\quad
W=\diag(\dfrac1{2b},\dfrac1{2a}),\quad D_0=G_M-W.
\]
The matrix $D_0$ is positive definite: its diagonal entries are
positive and $\det D_0=\dfrac1{4ab}$. With
$e_i'= \dfrac{F_i^+-1}{2a}+\dfrac{F_i^--1}{2b}\ge0$, the exact identity is
\[
\|x_i\|^2=\ell+z_i^TD_0z_i+e_i'.
\]
Completing the quadratic sums in the two groups gives
\begin{align}
\tr K&=F_1(R,c)+E_A+E_B,\label{eq:decomp1}\\
F_1(R,c)&=\dfrac{5R}{6}+\dfrac c3+\dfrac{25R}{12(1+c^2)},\notag
\end{align}
where $E_A,E_B\ge0$ are the sums of the nonnegative quadratic terms
and the terms $e_i'$ for $i=1,2$ and $i=3,4,5$, respectively. In particular
\[
E_B=\sum_{i=3}^5(z_i-z_*)^TD_0(z_i-z_*)+\sum_{i=3}^5e_i',
\quad z_*=(\dfrac{1}{3},-\dfrac{1}{3})^T.
\]
Let $A_0=R+2c=\|m_++m_-\|^2$, and $d_i=K_{ii}$. Put
\[
F_{12}=\det(\Z x_1+\Z x_2)^2=d_1d_2-K_{12}^2\ge1.
\]
Expansion of this coordinate-plane determinant gives
\[
d_1+d_2=\dfrac{A_0}{2}+\dfrac2{A_0}
 +\dfrac{2(F_{12}-1)}{A_0}+\dfrac{(d_1-d_2)^2}{2A_0}.
\]
Using this identity for $d_1+d_2$ gives
\begin{align}
\tr K&=F_2(R,c)+E_C+E_B,\label{eq:decomp2}\\
F_2(R,c)&=\dfrac{5R}{6}+\dfrac c3+\dfrac2{R+2c}
                           +\dfrac{4R}{3(1+c^2)},\notag
\end{align}
where $E_C\ge0$. The term $E_B$ is the same in both identities.

Next we bound $\tr K,a,b,c$.
Since $R\ge2\sqrt{1+c^2}$, equation \eqref{eq:decomp1} gives
\[
\tr K\ge f(c):=\dfrac{35+10c^2}{6\sqrt{1+c^2}}+\dfrac c3.
\]
On $[-\dfrac{1}{\sqrt3},0]$ its derivative
$\dfrac{1}{3}+\dfrac{c(10c^2-15)}{6(1+c^2)^{3/2}}$ is positive.
Moreover $f(-\dfrac{3}{10})=\dfrac{359}{6\sqrt{109}}-\dfrac{1}{10}>\dfrac{45}{8}$.
On $[0,\dfrac{1}{\sqrt3}]$, its second derivative has the sign of $4c^2-1$;
its derivative at the right endpoint is $-\dfrac{19}{48}$. Hence its minimum
on this interval occurs at an endpoint. Both $f(0)=\dfrac{35}{6}$ and
$f(\dfrac{1}{\sqrt3})=\dfrac{119\sqrt3}{36}$ exceed $\dfrac{45}{8}$. Therefore $c<-\dfrac{3}{10}$.

In fact $c>-\dfrac{1}{2}$. To see this put $x=-c$. If $x\ge\dfrac{1}{2}$, the
Lagrange reduction inequalities give $2x\le a\le\sqrt{1+x^2}$. The function
$a+\dfrac{1+x^2}{a}$ decreases on that interval, so
\[
A_0=R-2x\le\dfrac{x+\dfrac1x}{2}\le\dfrac{5}{4}.
\]
The last inequality is strict for $x>\dfrac{1}{2}$; at $x=\dfrac{1}{2}$ the first is
strict since $a>1$. This contradicts \eqref{eq:indicator} for the reference vector $u_A$.
Thus $\dfrac{3}{10}<x<\dfrac{1}{2}$.

Both $F_1$ and $F_2$ increase with $R$ on this range. For $F_2$,
\[
\partial_R F_2=\dfrac{5}{6}-\dfrac2{(R-2x)^2}+\dfrac4{3(1+x^2)}>\dfrac{23}{45}>0,
\]
using $R-2x\ge2(\sqrt{1+x^2}-x)>\dfrac{6}{5}$.
Put $t=\sqrt{1+x^2}$. Since $R\ge2t$, we have
\[
F_1(2t,-x)=\dfrac{5t}{3}-\dfrac x3+\dfrac{25}{6t},\quad
F_2(2t,-x)=\dfrac{8t}{3}+2\dfrac x3+\dfrac8{3t}.
\]
The first decreases with $x$ and the second increases, with derivative
$\dfrac{2}{3}+\dfrac{8x^3}{3t^3}>0$. They cross at
$x=\dfrac1{2\sqrt2}$, $t=\dfrac3{2\sqrt2}$, with common value
$\dfrac{71\sqrt2}{18}$. Equations \eqref{eq:decomp1}--\eqref{eq:decomp2}
prove
\begin{equation}\label{eq:smallbudget}
\tr K\ge \dfrac{71\sqrt2}{18},\quad 0\le E_B\le\dfrac{45}{8}-\dfrac{71\sqrt2}{18}<\dfrac{1}{20}.
\end{equation}
Also $F_2(2\sqrt{1+x^2},-x)$ at $x=\dfrac{5}{12}$ equals
$\dfrac{439}{78}=\dfrac{45}{8}+\dfrac{1}{312}$. Its monotonicity excludes $x\ge\dfrac{5}{12}$.
Consequently
\begin{equation}\label{eq:shape}
-\dfrac{5}{12}<c<-\dfrac{3}{10},\quad a,b<\dfrac{169}{144},\quad
A_0=R+2c<2+x^2-2x<\dfrac{149}{100}<\dfrac{3}{2}.
\end{equation}
Here $a+b<1+ab$ follows from $a,b>1$. This proves \eqref{eq:phaseinputs}.

We now estimate the determinant after orthogonal projection.
All sums in the remainder of this proof are over $i=3,4,5$.
Put $\eta_i=z_i-z_*$, so that $\sum\eta_i=0$, and set
\[
 F=\sum\eta_i^TD_0\eta_i,\quad E=\sum e_i',\quad
 \eps=F+E=E_B,\quad n_i=\eta_i^TW\eta_i,\quad N=\sum n_i.
\]
The two eigenvalues of
\[
 W^{-1/2}D_0W^{-1/2}
 =\begin{pmatrix}2t^2-1&-2xt\\-2xt&2t^2-1\end{pmatrix}
\]
are $(t-x)^2$ and $(t+x)^2$. Hence $N\le\dfrac{\eps}{(t-x)^2}$.
By the formulas for $F_i^\pm$ and $e_i'$, we also have
\[
 q_i=q_0+\ell_i-n_i+e_i',\quad
 q_0=\dfrac{8\ell}{9},\quad \ell_i=-2z_*^TW\eta_i,\quad \sum\ell_i=0.
\]
Direct inversion of the two-by-two matrix $D_0$ gives
\[
 (Wz_*)^TD_0^{-1}(Wz_*)=\dfrac{A_0-\ell}{9},\quad
 v:=\sum\ell_i^2\le\dfrac49(A_0-\ell)\eps.
\]
For a triple $w$ write
\[
 \mathcal D(w)=\dfrac34\Bigl(\sum w_i\Bigr)^2-\dfrac32\sum w_i^2.
\]
Since $y_3+y_4+y_5=0$, we have $3A_B^2=\mathcal D(q)$. Likewise $W^{1/2}\eta_3$,
$W^{1/2}\eta_4$, $W^{1/2}\eta_5$ sum to zero, so Heron's formula gives
$\mathcal D(n)\ge0$. Thus
\[
 \sum n_i^2\le \dfrac{N^2}{2},\quad
 \left|\sum\ell_i n_i\right|
 \le\sqrt{v\sum(n_i-\dfrac N3)^2}\le N\sqrt{\dfrac v6}.
\]
The derivative of $\mathcal D$ with respect to $w_i$ is
$\dfrac32(w_j+w_k-w_i)$. It is positive for every triple $w$ with
$w_i=q_0+\ell_i-n_i+\theta e_i'$, $0\le\theta\le1$. Indeed
\eqref{eq:smallbudget} gives $\eps<\dfrac{1}{20}$, and
\eqref{eq:shape} gives
\[
 q_0>\dfrac{4}{5},\quad N<\dfrac{9}{80}<\dfrac{1}{8},\quad E<\dfrac{1}{20},\quad
 v<\dfrac{13}{1000},\quad\max|\ell_i|<\dfrac{1}{10}.
\]
Here $\sum\ell_i=0$ gives $\max\ell_i^2\le\dfrac{2v}{3}$,
while $\ell\ge\dfrac1t>\dfrac{12}{13}$ and
$\dfrac49(A_0-\ell)<\dfrac49(\dfrac{149}{100}-\dfrac{12}{13})<\dfrac{13}{50}$.
Consequently $w_j+w_k-w_i>\dfrac{4}{5}-\dfrac{1}{5}-\dfrac{1}{8}-\dfrac{1}{20}>0$.
Dropping the terms $e_i'$ and expanding gives
\begin{equation}\label{eq:heron-quotient}
 3A_B^2\ge\dfrac94q_0^2-\dfrac32q_0N-\dfrac32v
                       -\sqrt{\dfrac{3}{2}}\,N\sqrt v.
\end{equation}
Here we used $\mathcal D(n)\ge0$.

The preceding bounds give
\[
 q_0\ge\dfrac8{9t},\quad
 N\le\dfrac{\eps}{(t-x)^2}\le\dfrac94\eps,\quad
 v\le\dfrac{13}{50}\eps.
\]
Here $(t-x)^{-2}=(t+x)^2<\dfrac94$, since
$t^2-x^2=1$, $x<\dfrac5{12}$, and $t<\dfrac{13}{12}$.
The right side of \eqref{eq:heron-quotient} increases with $q_0$
and decreases with $N,v$ on these ranges. Substitution therefore gives
\[
 3A_B^2\ge\dfrac{16}{9t^2}
 -\left[\dfrac3t+\dfrac{39}{100}
              +\dfrac94\sqrt{\dfrac{39\eps}{100}}\right]\eps.
\]
Since $t>1$ and $\eps<\dfrac1{20}$, the bracket is less than
\[
 3+\dfrac{39}{100}+\dfrac94\,\dfrac4{25}=\dfrac{15}{4};
\]
indeed $\dfrac{39}{2000}<\dfrac{16}{625}$.
Including $\eps=0$, we obtain
\[
 3A_B^2\ge\dfrac{16}{9t^2}-\dfrac{15}{4}\eps.
\]

By \eqref{eq:decomp1}--\eqref{eq:decomp2} and $\eps=E_B$, we have
\[
 \tr K-\eps\ge\max\{F_1(R,c),F_2(R,c)\}\ge\dfrac{71\sqrt2}{18}>\dfrac{251}{45}.
\]
The last comparison follows by squaring, with difference $\dfrac{23}{4050}$.
Put $s=\dfrac{45}{8}-\tr K>0$.
Then $\eps<\dfrac{17}{360}-s$, and $t<\dfrac{13}{12}$ gives
\[
 2A_B^2-\dfrac89
 \ge\dfrac{32}{27t^2}-\dfrac89-\dfrac52\eps
 >\dfrac{184}{1521}-\dfrac{17}{144}+\dfrac52s
 =\dfrac{71}{24336}+\dfrac52s>0.
\]
This proves \eqref{eq:d5-area-gap}, and therefore
Proposition~\ref{clm:d5-phase}.

\end{proof}

\begin{proof}[Proof of Proposition~\ref{clm:d5-dual}]
For $r\in\mathcal R$, write
\[
R_r=\|K^{-1/2}r\|^2=r^TK^{-1}r.
\]
Put $r_-=e_1-e_2$ and $R_-=R_{r_-}$, so that $c_-=c_{r_-}$.
From \eqref{eq:floors} and $\det\Lambda=\dfrac{\sqrt{\det K}}{2}$,
$\lambda_1(\Lambda^*)\ge\dfrac{4}{3\sqrt{\det K}}$. Hence every root satisfies
\begin{equation}\label{eq:allroots}
R_r\ge\lambda_1(\Lambda^*)^2\ge\dfrac{16}{9\det K}.
\end{equation}
Taking the trace of the root expansion, and of its product with
$K^{-1}$, gives
\begin{equation}\label{eq:weightedtrace}
2\sum_r c_r=\tr K,\quad \sum_r c_rR_r=5.
\end{equation}

Only $e_1+e_2$ and $e_1-e_2$ contribute to $K_{12}$.
Thus
$2c_-\ge-2K_{12}=(K_{11}+K_{22})-A_0$.
Note that
$F_{12}=K_{11}K_{22}-K_{12}^2\ge1$. Since
$K_{12}=\dfrac{A_0-K_{11}-K_{22}}{2}$ and $K_{11}K_{22}\le\dfrac{(K_{11}+K_{22})^2}{4}$,
\[
1\le\dfrac{A_0(K_{11}+K_{22})}{2}-\dfrac{A_0^2}{4},\quad
(K_{11}+K_{22})\ge \dfrac{A_0}{2}+\dfrac2{A_0}.
\]
Proposition~\ref{clm:d5-phase} gives $A_0<\dfrac{3}{2}$, so
\[
2c_-\ge\dfrac2{A_0}-\dfrac{A_0}{2}>\dfrac{7}{12}.
\]

The reference vectors $(h_+,h_-,e_1,e_3,e_4)$ form a basis of
$\Lambda_0^*$: one recovers
$e_2=h_++h_--e_1$ and $e_5=h_+-h_--e_3-e_4$, and also
$\dfrac{e_1+\cdots+e_5}{2}=h_+$.
Thus $(m_+,m_-,x_1,x_3,x_4)$ is a basis of $\Lambda$.
The primitive dual vector $K^{-1/2}r_-$ pairs with this basis as
$(0,0,1,0,0)^T$. Its kernel section has basis
$(m_+,m_-,x_3,x_4)$.
Projecting this section modulo $M$ gives squared determinant $A_B^2$,
since $\det M=1$. The primitive dual-section formula used in
\eqref{eq:short-dual} gives the same squared determinant as
$\det(\Lambda)^2\|K^{-1/2}r_-\|^2=(\dfrac{\det K}{4})R_-$.
Therefore
\begin{equation}\label{eq:dual-section-area}
(\det K)R_-=4A_B^2.
\end{equation}

Separating the term $r=r_-$ in \eqref{eq:weightedtrace}, we obtain
\[
\begin{aligned}
5\det K
&=(\det K)\sum_r c_rR_r\\
&\ge\dfrac{16}{9}\sum_r c_r+
 c_-\left((\det K)R_--\dfrac{16}{9}\right)\\
&=\dfrac89\tr K+2c_-\left(2A_B^2-\dfrac89\right).
\end{aligned}
\]

This proves \eqref{eq:new-d5-trace} and completes the proof.
\end{proof}

\paragraph{\textbf{AI usage. }}ChatGPT assisted with the calculation and the composite in the $D_5$ case in the appendices. 

\begingroup\raggedright

\endgroup
\end{document}